\documentstyle[11pt,leqno]{article}
\input boxedeps.tex
\SetepsfEPSFSpecial
\HideDisplacementBoxes
\def\reals{\hbox{\sl I\kern-.18em R \kern-.3em}}
\newtheorem{theorem}{Theorem}
\newtheorem{lemma}{Lemma}
\def\np{\noindent}
\def\pf {\np {\bf Proof:} \ }
\def\endpf{$\|$ \bigskip}
\def\cA{{\cal A}}
\def\cH{{\cal H}}
\def\cU{{\cal U}}
\def\cX{{\cal X}}
\def\cV{{\cal V}}
\def\cW{{\cal W}}
\def\pp{{\prime\prime}}

\begin{document}

\title{On Markov's Theorem}

\author{Joan S. Birman%
\thanks{The first author acknowledges partial
support from the U.S.National Science
Foundation, under Grant DMS-9973232.}
\\e-mail jb@math.columbia.edu \and William W. Menasco
\thanks{The second author acknowledges partial
support from the U.S. National
Science Foundation under grant DMS-9626884}
\\e-mail menasco@tait.math.buffalo.edu}
\date{February 15, 2002 }
\maketitle
\centerline {to appear in Proceedings KNOTS-2000} 
\centerline {(special issue, Journal of Knot Theory and its
Ramifications)}

\section{Introduction}
\label{section:introduction}
Let $\cX$ be an oriented link type in
the oriented 3-sphere $S^3$ or $\reals^3 = S^3 - \{\infty\}.$  A
representative $X \in \cX$ is said to be a {\em closed braid} if
there is an unknotted curve \glossary{}{${\bf A}
\subset S^3 - X$} (the {\em axis}) and a choice of fibration $\cH$ 
of the open solid
torus $S^3 - {\bf A}$  by meridian discs
$ \{ H_{\theta}: \ \theta \in [0,2\pi] \}$,  such that whenever
$X$ meets a fiber
$H_\theta$ the intersection is transverse.
The fact that the link
$X$ is a closed braid implies that the number of points in
$X \cap H_\theta$ is independent of
$\theta$, and we call this number the {\em braid index} of $X$. The braid index
of $\cX$ is the minimum value of the braid index of $X$
over all closed braid representatives $X \in \cX$.

\begin{figure}[htpb]
\centerline{\BoxedEPSF{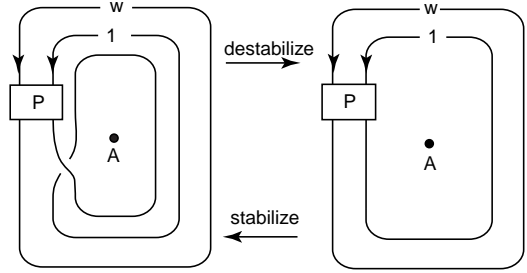 scaled 600}}
\caption{Stabilization and destabilization of closed braids.}
\label{figure:stab-destab}
\end{figure}

Closed braid representations of $\cX$ are not unique, and Markov's well-known 
theorem \cite{Markov} asserts that any two are related by a finite
sequence of elementary moves. One of them is {\em braid isotopy}, by which we
mean isotopy in the complement of the braid axis which preserves transversality
between $X$ and fibers of
$\cH$.   The
other two moves are mutually inverse, and are illustrated in Figure
\ref{figure:stab-destab}. Both take closed braids to closed braids. We call them
{\em stabilization} and {\em destabilization}, where the former increases the
braid index by one and the latter decreases it by one.
The `weight' $w$ denotes $w$ parallel strands, relative to the given
projection. The braid inside the box which is labeled $P$ is
an arbitrary $(w+1)$-braid.  

\begin{theorem}[Markov's Theorem:]
\label{theorem:MT}
Let $X_1,X_2$ be closed braid representatives of the same oriented link type
$\cX$ in oriented 3-space, with the same braid axis ${\bf A}$.  Then
$X_2$ may be obtained from $X_1$ by braid isotopy and a finite number of
stabilizations and destabilizations.
\end{theorem}

We know of 4 published proofs of Theorem \ref{theorem:MT}, see
\cite{Birman1974,Lambropoulou-Rourke,Morton1986b,Traczyk}, and each gives a new
way of looking at the result. The main result in this paper is to give yet
another proof!  We hope that our new proof will be of interest because it gives
new insight into the geometry, as follows.  During recent
years a fairly clear picture has emerged of the closed braid representatives of
the unknot and unlink. See the main theorem in \cite{BM-5}.
Since the operation of taking the braid connected sum of any closed braid
representative $X$ of $\cX$ with any closed braid representative $U$ of the unknot
produces a different closed braid representative of $\cX$, a
natural question is whether this process of taking braid connected
sums with copies of closed braid representatives of $\cU$ explains all of the
complications in closed braid representatives of $\cX$?  Our proof of
Theorem \ref{theorem:MT} will clarify this situation.  Let $X,X'$ be any two
closed braid representatives of the same knot or link type $\cX$. We will show
that there is an isotopy taking $X_1$ to an intermediate closed braid
$X_3$ and another isotopy from $X_3$ to
$X_2$  such that:
\begin{enumerate}
\item  $X_3$ is obtained from $X_1$ by taking the braid
connected sum of $X_1$ with some number of copies of closed braid representatives
of $\cU$.  
\item The isotopy that takes $X_3$ to $X_2$ is a
push across an embedded annulus $\cA_2$ which is a subset of a Seifert
surface ${\bf F}$ for $X_3$. In particular $X_2$ is a
preferred longitude for $X_3$. 
\end{enumerate}
The proof of Theorem \ref{theorem:MT} is based upon ideas in a forthcoming paper by
the authors \cite{BM-8}, in which `Markov's Theorem Without Stabilization' will be
proved. In that paper we will describe a set of moves which suffice to take any
closed braid representative of a knot or link to one of minimum braid index 
for the link type $\cX$, with each move preserving or reducing the braid index.  The
proof which is given here of Theorem \ref{theorem:MT} will also be used by Nancy
Wrinkle in her PhD thesis, which concerns knots which are transverse to the standard
contact structure in $\reals^3$.   It is our hope that the geometry revealed in the
new proof of Theorem \ref{theorem:MT} will have other applications too.

\section{Braid foliations}
\label{section:braid foliations}
The main result in this paper (Theorem
\ref{theorem:MT}) is about the relationship between two closed braid
diagrams which represent the same link. However the work which we will do to
prove it will not be in the setting of knot diagrams. Rather, we will be dealing
with surfaces which our links bound, and with certain {\it braid foliations} of
these surfaces. The foliated surfaces have been used before, in our 6 earlier
papers with the common title ``Studying links via closed
braids", e.g.see \cite{BM-5}. In this section we will review and describe
the machinery which  we need from these other papers.  We will refer the reader to the
review article \cite{B-F} for proofs, whenever possible.  

We are given a representative $X$ of an oriented link type $\cX$, where $X$ is a
closed $n$-braid with braid axis ${\bf A}$.  Let $\cH =
\{H_{\theta}; 0 \leq \theta \leq 2\pi \}$ be a choice of disc fibers of the braid axis
complement, where $H_{\theta}$ denotes a fiber of $\cH$. The braid axis and the
fibers of $\cH$ will
be seen to serve the role of a coordinate system in 3-space. We will use it to
describe the link
(and an auxiliary surface which it bounds) via a set of combinatorial data. The
fact that
$X$ is a closed braid with respect to $\cH$ implies that it intersects each fiber
$H_{\theta}$ of $\cH$
transversally in exactly $n$ points. The closed braid $X$ is oriented so that it is
pointing in the
direction of increasing $\theta$ at each point of $X \cap H_\theta$. 
\begin{figure}[htpb]
\centerline{\BoxedEPSF{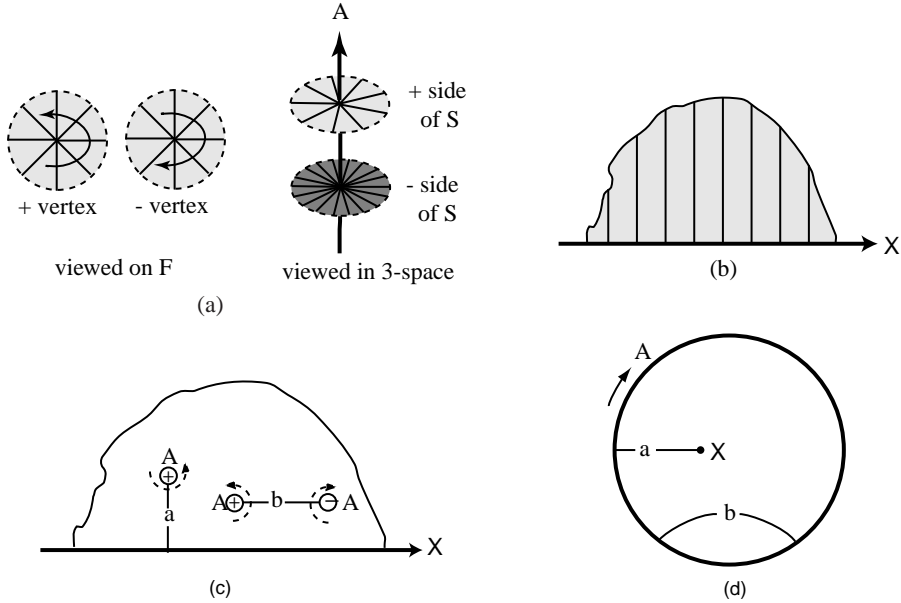  scaled 750}}
\caption{(a)Foliation of ${\bf F}$ near a vertex; (b) Foliation of ${\bf F}$ near
$X=\partial{\bf F};$ (c) $a$-arcs and $b$-arcs, viewed on ${\bf F}$;  (d) $a$-arcs and
$b$-arcs, viewed on a fiber
$H_\theta$.}
\label{figure:foliations1}
\end{figure}

Our link is assumed to be the boundary of a surface \glossary{}{${{\bf F}}$} of 
maximum Euler characteristic. After
modifying ${{\bf F}}$ we may assume that it admits a special type of singular foliation
which was studied and used by the authors in their earlier papers and reviewed in
\cite{B-F}. The foliation is radial in a neighborhood of each point of ${\bf A} \cap
{{\bf F}}$ (see Figure
\ref{figure:foliations1}(a)) and
transverse to the boundary in a neighborhood of $\partial {{\bf F}}$ (Figure 
\ref{figure:foliations1}(b)).  Notice that in Figure
\ref{figure:foliations1} the surface ${{\bf F}}$ is naturally oriented
by the orientation on $X$, which is
chosen so that the polar angle $\theta$ is strictly increasing as we walk along
$X$.  This fact orients the associated flow and could be used to orient the
foliation, although we have not done so. A {\it vertex} in the foliation is a point
in
${\bf A}\cap{\bf F}$.  The braid axis ${\bf A}$ pierces ${{\bf F}}$ from either the
negative or the positive side at each pierce
point, and we have indicated this by attaching positive or negative signs to the
 pierce points on
${{\bf F}}$. A {\em leaf} in
the foliation is a component of  intersection of $H_{\theta}$ with the surface 
${{\bf F}}$.  We may
assume that each leaf is an arc. For details on this assertion and others like it see
\cite{B-F}.
Leaves are {\it singular} if they contain a
singularity of the foliation, otherwise they are {\it non-singular}.
The singularities may be
assumed to be finite in number and to occur on distinct
fibers of $\cH$.  Every singularity may be
assumed to result from a saddle point tangency between ${{\bf F}}$
and a fiber of $\cH$. Non-singular
leaves are said to be {\it type a} (respectively {\it type b}) if they are arcs
which have one
endpoint on ${\bf A} \cap {{\bf F}}$ and the other on $X = \partial {{\bf F}}$
(respectively both endpoints
on ${\bf A} \cap {{\bf F}}$). See Figure
\ref{figure:foliations1}(c) and (d). Leaves which have both their endpoints on
$X$ do not occur because $ {\bf F} $ is orientable (see Lemma 1.1 of \cite{B-F}). 
When the foliation of ${\bf F}$ has all of these properties we call ${\bf F}$ a {\em
Markov surface}. 

If a vertex is the endpoint of an $a$-arc, it is always positive
(as indicated in
Figure \ref{figure:foliations1}(a),
but it could be either positive or negative at the endpoint of
a $b$-arc. To indicate these differences we will sometimes show a point where the
axis pierces ${\bf F}$ as as a circle with a $\pm$ sign inside it. 
Another way to think of the signs on the vertices is that they
are {\it positive} or {\it negative} according as the outward-drawn oriented normal to 
${\bf F}$ has the same
or opposite orientation as the braid axis at the vertex. This means that when we
 view the positive
side of
${{\bf F}}$, the sense of increasing $\theta$ around a vertex will be counterclockwise
(resp.
clockwise) when the vertex is positive (resp. negative). The {\it valence} of a
 vertex is the number of singular leaves which have endpoints at that vertex.

\begin{figure}[htpb]
\centerline{\BoxedEPSF{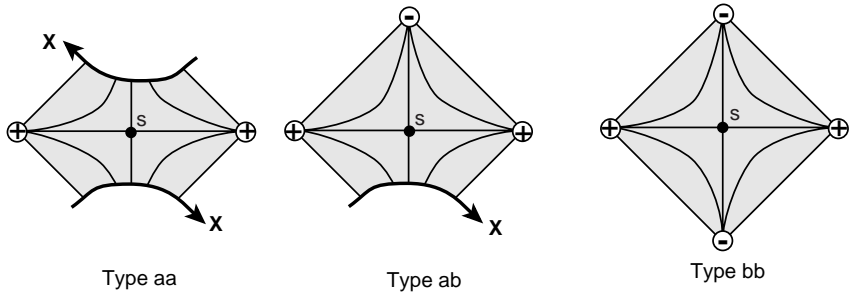 scaled 650}}
\caption{The three tile types.}
\label{figure:foliation of tile types}
\end{figure}

The foliation may be used to decompose the surface ${{\bf F}}$ into
a union of foliated $2$-cells, each of which contains exactly one singularity of
 the foliation.  We refer to these 
$2$-cells as {\it tiles} and the resulting decomposition of ${{\bf F}}$ as a {\it
tiling}. Each tile is a regular neighborhood on
${\bf F}$ of its singular leaves.  See Figure \ref{figure:foliation of tile types}. 

The {\it
vertices} of the tiles are the points where the braid axis ${\bf A}$ intersects
the  surface ${{\bf F}}$.
(There are also other vertices on the graph of singular leaves, but we prefer
to exclude them when we refer to tile vertices.) The {\it edges} of the 
tiles are non-singular leaves of type $a$ or $b$. (There are also other tile edges
which are  subarcs of $X$,
but it will be convenient to ignore those too, just as we ignored the vertices 
which are on $X$.)

If two tiles intersect, then they intersect along a non-singular leaf of
the foliation, which is then necessarily an edge of type $a$ or $b$. The tiles 
fall into three
types, according to their  foliations. We call them types
$aa, ab$ and $bb$, the notation indicating that in an $aa$-tile (respectively 
$ab, bb$-tile) the
leaves which come together before the singularity are both type $a$ 
(respectively types $a$ and
$b$, both type $b$).
The three tile types are distinguished by the number of their vertices, the types
 of their
non-singular leaves, and the endpoints of their singular leaves. Tiles of type $ab$
 and $aa$ meet
$X=\partial{\bf F}$, but tiles of type $bb$ are contained entirely in the interior
of ${\bf F}$.  Notice
also that tiles of type $bb, ab, aa$ each contain $2$ positive vertices, but
contain $2,1$ and $0$
negative vertices respectively. The segments where $X$ meets a tile of type $ab$
 or $aa$ will be seen to `act like a negative vertex'.

The singularities also have signs. Let $s$ be a
singular point of the foliation of a Markov surface ${{\bf F}}$ for a link $X$, and
let
$H_\theta$ be the disc fiber which contains $s$.  We say that $s$ is {\it
positive} if the outward-drawn oriented normal to the
oriented surface ${{\bf F}}$ coincides in direction with the normal to
$H_\theta$ in the direction of increasing $\theta$. Otherwise $s$ is {\it
negative}.  

\begin{figure}[htpb]
\centerline{\BoxedEPSF{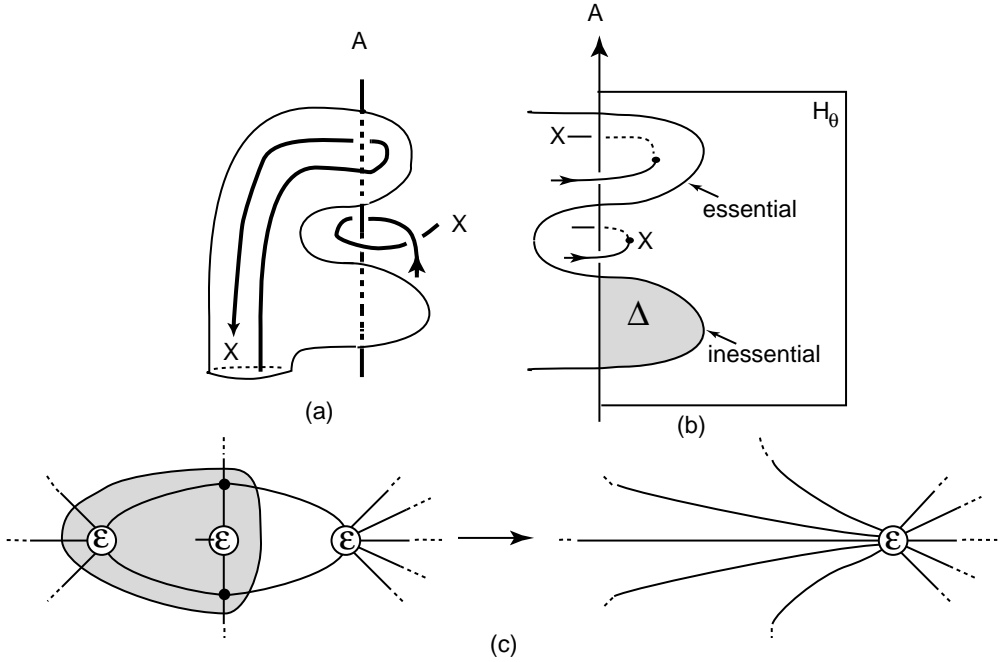 scaled 750}}
\caption{Essential and inessential b-arcs. The bottom two sketches show the changes
in typical singular leaves in the foliation when an inessential $b$-arc is removed.}
\label{figure:essential b-arcs}
\end{figure}

A $b$-arc, when viewed on a fiber of $\cH$, divides the fiber into
two discs.  We say that
it is {\em essential} if both discs are pierced by $X$.  This condition can be seen 
in the ordering along ${\bf A}$ of the vertices in the tiling of ${\bf F}$. If the
two vertex endpoints of the $b$-arc have adjacent numbers in the cyclic order on
${\bf A}$, then the $b$-arc is inessential, otherwise it is essential. See Figure
\ref{figure:essential b-arcs} (a) and (b).
\begin{lemma}
\label{lemma:b-arcs essential}
All $b$-arcs may be assumed to be essential.  
\end{lemma} 

\pf Any inessential $b$-arc may be removed by an isotopy of ${\bf F}$. For details see
\cite{B-F}.  The isotopy removes two vertices and two singularities from the
foliation. See Figure \ref{figure:essential b-arcs} (c). \endpf

We next describe two moves which
modify the closed braid $X$ and the foliation on ${\bf F}$, but do not change the link
type $\cX$.  Both  change the braid index.  Stabilization adds 
a trivial loop and destabilization removes it.
However, there is more to it then that. Our stabilization and
destabilization moves are guided
by the foliation of the surface. Thus they are moves of the pair
$(X,{\bf F})$. As a move on the
pair $(X,{\bf F})$ our stabilization move is {\em not} the
inverse of our destabilization move, even
though they {\em are} mutually inverse in terms of their effect on $X$.
\begin{figure}[htpb]
\centerline{\BoxedEPSF{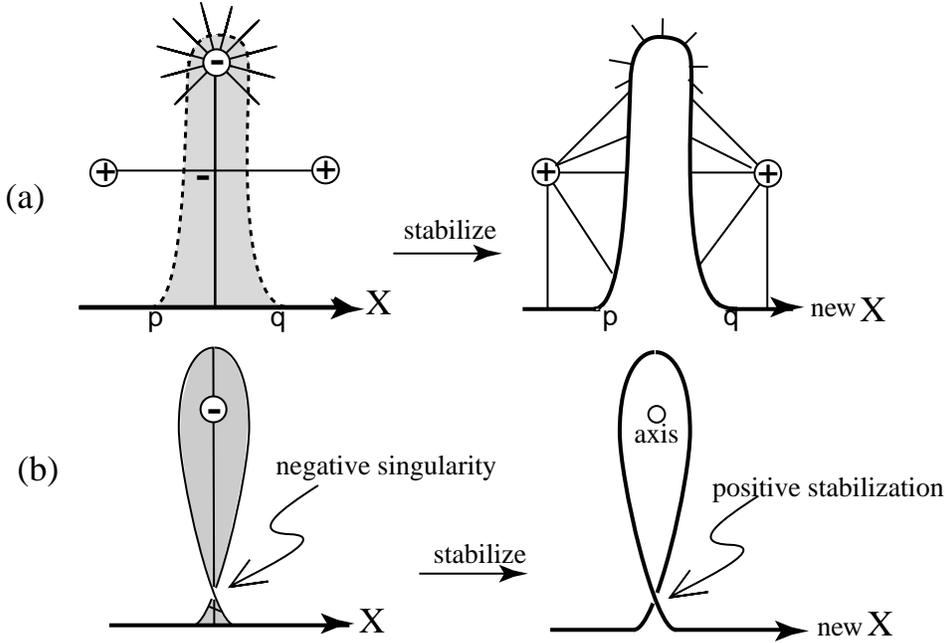 scaled 750}}
\caption{Stabilization along an ab-tile.}
\label{figure:stabilization along ab-tile}
\end{figure}

\

\np{\bf Stabilizing along ab tiles:}
Figure \ref{figure:stabilization along ab-tile}
shows a modification in the tiling on ${{\bf F}}$, which we call
{\em stabilizing along an $ab$ tile}. It is realized
by pushing $X$ along a disc neighborhood
of the singular leaf in an $ab$ tile. The boundary of the neighborhood may be 
chosen to be
transverse to the foliation, so the result is a new closed braid representative
of the $X$,
bounding a new surface ${{\bf F}}^{\prime}$ which is tiled, but with one less tile 
than ${{\bf F}}$.
We may visualize the result of the isotopy of $X$ as the
addition of a trivial loop, as depicted in Figure \ref{figure:stabilization along ab-tile}.
This move increases the number of braid
strands. As illustrated in Figure \ref{figure:stabilization along ab-tile},
its effect on the braid representation
of $X$ is to add a
``trivial loop'' around the axis, increasing the braid index by 1. The effect of
our move on the
tiling is to eliminate a negative vertex. This will change any $bb$ tile that is
adjacent to the negative vertex to an $ab$ tile.

\

\begin{figure}[htpb]
\centerline{\BoxedEPSF{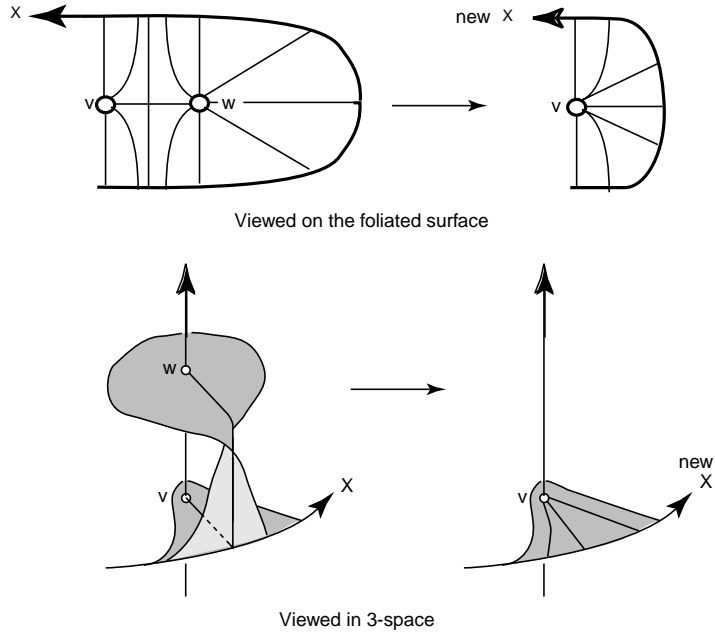 scaled 600}}
\caption{Destabilization along an end tile.}
\label{figure:destabilization along end tile}
\end{figure}

\np {\bf Destabilizing along end tiles:}
If the tiling has a vertex of
valence 1, then that vertex must be in an $aa$-tile which is glued to itself 
(see Figure \ref{figure:destabilization along end tile}),
forming a trivial loop.  Call such a tile an {\em end tile}. Removing it,
we obtain a new closed braid, with braid index one less than that of the original
one. We call
this process destabilizing along an end tile. This move is the inverse of
stabilization along an $ab$ tile, as far as its effect on the closed braid, but
the two moves are not mutually inverse as regards their effect on the foliation.

\section{A new proof of Markov's Theorem}
\label{section:Markov's theorem}
In our proof the  closed braid $X$ will be assumed to be the boundary of a Markov
surface ${\bf F}$. The central object of study is the pair $(X,{\bf F})$ 
rather than simply the closed braid $X$. The reader is referred to 
$\S$\ref{section:braid foliations} for the
description of the stabilizing and destabilizing moves, as they 
are reflected in changes
in ${\bf F}$. They play a central role in our proof.

The new proof of the MT was suggested to us by ideas which were sketched by 
D. Bennequin in \cite{Bennequin}, however his proof made use of the 
characteristic foliation of a surface bounded by a transverse knot, rather than
the braid foliations which were described in this paper. We read his proof of
the MT many times before we wrote this paper but we only arrived at our present
understanding of it comparatively recently, because of its very  sketchy nature.
(We may, however, have misunderstood what he intended to say.) 

If $X$ and $X'$ are closed
braids which represent the same link type $\cX$, we say that 
$X$ and $X'$
are {\em Markov-equivalent} (written $X\equiv X'$) if $X$ can be changed to 
$X'$ by braid isotopies together with a finite number of 
stabilizations and destabilizations, i.e. the allowed moves in
the MT. 

\begin{lemma}
\label{lemma:MT is true for the unknot}
Let \glossary{}{$U$} be an arbitrary closed braid representative of 
the unknot
and let \glossary{}{$U_0$} be the standard 1-braid
representative. Then $U\equiv U_0$.  \end{lemma}

\pf We are given a closed braid representative $U$ of the
unknot. It is the boundary of a disc $D$, and by the methods 
which are described in
$\S$\ref{lemma:b-arcs essential}
we may assume that the disc has a tiling, and that all 
$b$-arcs in the tiling are
essential. There will in general be $p$ positive vertices and 
\glossary{}{$n$} negative vertices in
the tiling. Notice that the braid index is $p-n$, because the 
braid index is the linking
number of $U$ with the braid axis ${\bf A}$, but this linking number 
may also be computed as
the algebraic intersection number of ${\bf A}$ with a disc $D$
which $U$ bounds, i.e.
$p-n$. 

If $n>0$ then there exist $ab$ and/or $bb$ tiles, because there are 
no negative
vertices on an $aa$ tile.  However a $bb$ tile can only be adjoined 
to another $bb$
tile or an $ab$ tile, but since $bb$ tiles do not meet the boundary 
it is impossible
to have only $bb$ tiles. So if $n>0$, then there exists an $ab$ tile. 
Let us stabilize
along it, as described in $\S$\ref{section:braid foliations}
and in Figure \ref{figure:stabilization along ab-tile}. This 
eliminates one negative
vertex (at the expense of increasing the braid index). 
The process ends when there are
no more negative vertices, i.e. when $D$ is tiled entirely with $aa$
tiles.

Assuming ${\bf D}$ is tiled entirely by $aa$ tiles, we consider the graph $G$ 
which is the union of
the non-separating edges in the singular leaves and the  vertices in the
tiling. The graph \glossary{}{$G$} must be a tree, because $G$ is a deformation 
retract of $D$ and
$D$ is a disc. Choose an end-point vertex of this tree and 
eliminate it by
destabilizing, as in Figure \ref{figure:destabilization along end
tile}. This gives a
new disc which is again tiled entirely by $aa$ tiles, but with one 
less tile and so one
less singularity. The process ends with a subdisc $D$, which 
is foliated
radially without singularities. Its boundary is the standard 
1-braid representative 
$U_0$ of the unknot. \endpf

Lemma \ref{lemma:MT is true for the unknot} is a very special cases of the MT. Two
other special cases follow easily:    

\begin{lemma}
\label{lemma:MT true for braid connected sum with unknot}
Let $X$ and $U$ be links which are presented as closed
braids. Assume that $U$ represents the unknot. Then $X$ is 
Markov-equivalent to the
braid-connected sum of $X$ and $U$.
\end{lemma}

\pf The closed braids $X$ and $U$ bound Markov surfaces
${\bf F}$ and $D$, where $D$ is a disc. The surface ${\bf F} \#
D$  is then a Markov surface whose boundary is the braid connected sum of $X$ and
$U$.  We may simplify the tiling on $D$ by stabilizing along $ab$-tiles and then
deleting  valence one
vertices, as in the proof of Lemma \ref{lemma:MT is true for the unknot}. 
This produces a Markov equivalence between $X \# U$ and $X$.   \endpf

\begin{lemma}
\label{lemma:MT true for preferred longitude}
Let $X$ and $X'$ be closed braids, in general
having distinct braid indices. Assume
that $X'$ is a preferred longitude for $X$.  Then
$X\equiv X'$. 
\end{lemma}

\pf Since $X'$ is a preferred longitude for $X$ there is a Seifert surface
${\bf F}$ with $\partial{\bf F} = X$  and an annulus
$\cA\subset{\bf F}$, with $\partial\cA = X - X'$. 
(Here the minus sign denotes a reversal of the
order of $X'$.) Using the tools which were reviewed in
$\S$\ref{section:braid foliations} above, we may assume that
${\bf F}$ has maximal Euler characteristic and  is a Markov surface. Therefore
${\bf F}$ admits a tiling.  Since $X'$ is a closed braid we may further assume that
$X'\subset{\bf F}$ is everywhere transverse to the leaves in the foliation. Our
plan is to use the tiling on ${\bf F}$ and the induced tiling on $\cA\cap{\bf F}$
to prove that  $X\equiv X'$.  

\begin{figure}[htpb]
\centerline{\BoxedEPSF{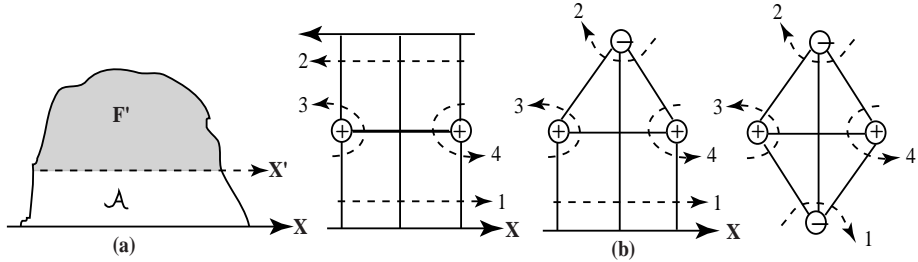 scaled 700}}
\caption{Passages of $X'$ through (a) the annulus bounded by $X$ and $X'$
and (b)tiles of type $aa$, $ab$ and $bb$ in the annulus.}
\label{figure:tiles in the annulus}
\end{figure}

Recall that ${\bf F}$ is a union of tiles of type $aa, ab$ and $bb$. 
Let us examine how $X'$
intersects the tiles in the tiling of
${\bf F}$. Since ${\bf F}$ is tangent to fibers of $\cH$ at its singular points, 
whereas $X'$
(being a closed braid) is everywhere transverse to fibers of $\cH$, 
it follows that
$X'$ never passes through a singular point in the tiling, also 
it is oriented so
that it travels clockwise (resp. counterclockwise) whenever it 
goes around a negative
(resp. positive) vertex in the tiling, when we are viewing the 
positive side of ${\bf F}$.
Also, whenever subarcs of $X'$ are locally parallel to subarcs of 
$X$ these subarcs are
oriented coherently. It follows that for each of the three tile types 
there are four
different ways that $X'$ can pass through a tile $T$. These are illustrated 
in Figure \ref{figure:tiles in the annulus}(b).
We have assigned little numbers $1,2,3,4$ to the four possible 
components of  $X'
\cap T$ in each tile $T$.

\begin{figure}[htpb!]
\centerline{\BoxedEPSF{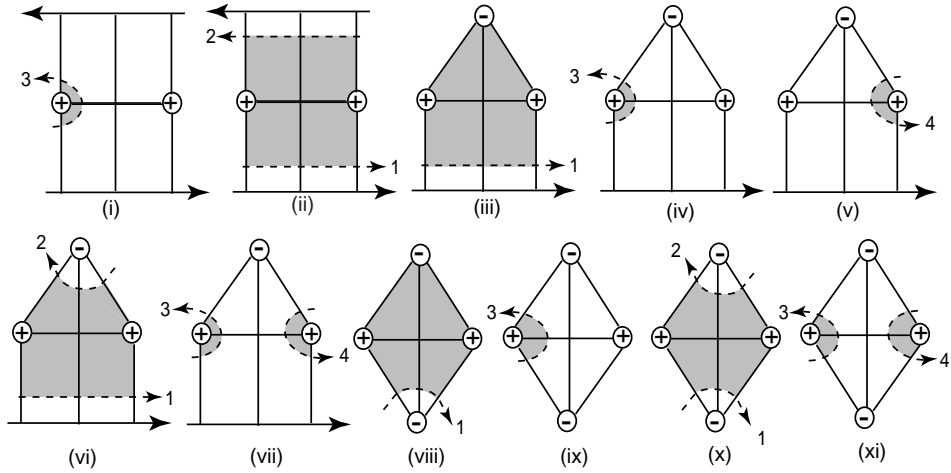 scaled 650}}
\caption{Passages of $X'$ through tiles in the annulus: 11 cases.}
\label{figure:tiles in annulus}
\end{figure}

The surface ${\bf F} = {\bf F}^\prime \cup \cA$ is oriented, also $X'$ is
oriented so that ${\bf F}^\prime$ is always on its left and $\cA$ on its right,
as in Figure \ref{figure:tiles in the annulus},
whereas $X$ always has $\cA$ on its left.
A check of the possibilities then shows that if $X' \cap T$ has more than one
component, then the components can only be arcs of type 1 and 2 or type 3 and 4. In
the case of tiles of type $aa$ there is an additional constraint which implies that
we cannot have components of type 3 and 4 simultaneously. For, if this occurred,
then one of the two singular arcs in the tile would begin and end on $X$ and
the other would begin and end on $X'$ and both arcs would be in $\cA$. However
this is impossible because these arcs intersect once in the tile,
however $\cA$ is an annulus and arcs of the type we have just described must
intersect an even number of times. Thus after removing obvious duplications, we are
reduced to the 11 cases which are illustrated in Figure 
\ref{figure:tiles in annulus}.

One possible measure of the complexity of $\cA$ is the number of negative
vertices in $\cA$.  If a tile $T$ has a negative vertex and intersects
$X = \partial{\bf F}$ then $T$ must be type $ab$.  Let us
travel along $X$, looking for $ab$ tiles. If one occurs then one of its 
singular
leaves begins at $X$ ends at a negative vertex. We may eliminate that 
negative
vertex by stabilizing along the $ab$ tile (see Figure
\ref{figure:stabilization along ab-tile}). 
This
increases the braid index by 1. The sign of the singularity in the $ab$ tile will
determine the sign of the half-twist which we add.
This process ends when there are no
more singular leaves which have one endpoint on $X$ and the other on a 
negative vertex.
But then, cases (iv),(v),(vii) do not occur. In fact, 
cases (vi),(viii),(ix),(x) and (xi) also do not 
occur because each of those presupposes the 
existence of a negative vertex in
$\cA$, however if we begin at a negative vertex in $\cA$ and 
travel along singular leaves
without leaving $\cA$ we must eventually arrive at a singular 
leaf in $\cA$ which has one
endpoint at a negative vertex and the other on $X$, but that is impossible.   

The only remaining cases are cases (i),(ii),(iii). If case (i) occurs we may
stabilize $X'$ by pushing it across the singular leaf in ${\cA
}$. (Equivalently, we could destabilize $X$ by repeatedly deleting vertices of
valence 1). This pushes a positive vertex from $\cA$ to ${{\bf F}^\prime}$. But
then, after some number of such stabilizations, we will have 
reduced to the situation
where only cases (ii) and (iii) occur. That is, there are no singularities
whatsoever in $\cA$, so the modified braid representatives of $X$ and $X'$
have identical braid structures. The proof of 
Lemma \ref{lemma:MT true for preferred longitude} is complete. \endpf

Lemma \ref{lemma:MT true for preferred longitude} 
was not the general case of the MT because we assumed that $X'$ is a
preferred longitude for $X$.  The next lemma will be needed
before we can attack the general case. 

\begin{lemma}
\label{lemma:unlinking}
{\bf The unlinking lemma \cite{Bennequin}.} Let $X$ and $X^{\prime}$ be links in
3-space $\reals^3$ which are not separated by any topological plane
$\reals^2\subset\reals^3$. Then we may modify the link
$X^\prime$ by  isotoping it to a link $X^\pp$
which is separated from $X$ by a  plane $P$ in $\reals^3$.
The isotopy from $X^\prime$ to $X^{\prime \prime}$ may 
be realized by choosing disjoint embedded discs $D_1,\dots,D_k$ whose interiors are
disjoint from $X'$ and $X^{\prime \prime}$, with $int(D_i)\cap X$ a point. Each
$D_i$ intersects $X$ in an arc $\alpha_i \subset \partial D_i$ and
 $X^{\prime \prime}$ in the closure $\beta_i$ of the complement of
$\alpha_i$ in $\partial D_i$. The isotopy is a push of 
$\alpha_1,\dots,\alpha_k \subset X^\prime$ across 
$D_1,\dots,D_k$ to $\beta_1,\dots,\beta_k \subset X^{\prime \prime}$. 
\end{lemma}

\pf We may assume that $X\cup X'$ is in general
position with respect to projection onto a plane which is orthogonal
to the braid axis. Think of $X$ as being colored red and $X'$ as being colored
green. There are finitely many points where the projected image of green
crosses that of red, with each crossing transversal. If, at  {\it every} such
crossing, green crosses under red, then
$X$ and $X'$ will be unlinked geometrically and there is no obstruction to pushing
green below red to separate them. If not, then after a finite number of
green-red crossing switches they will be unlinked. Each crossing
switch may clearly be described as in the statement of the lemma. \endpf

\begin{figure}[htpb]
\centerline{\BoxedEPSF{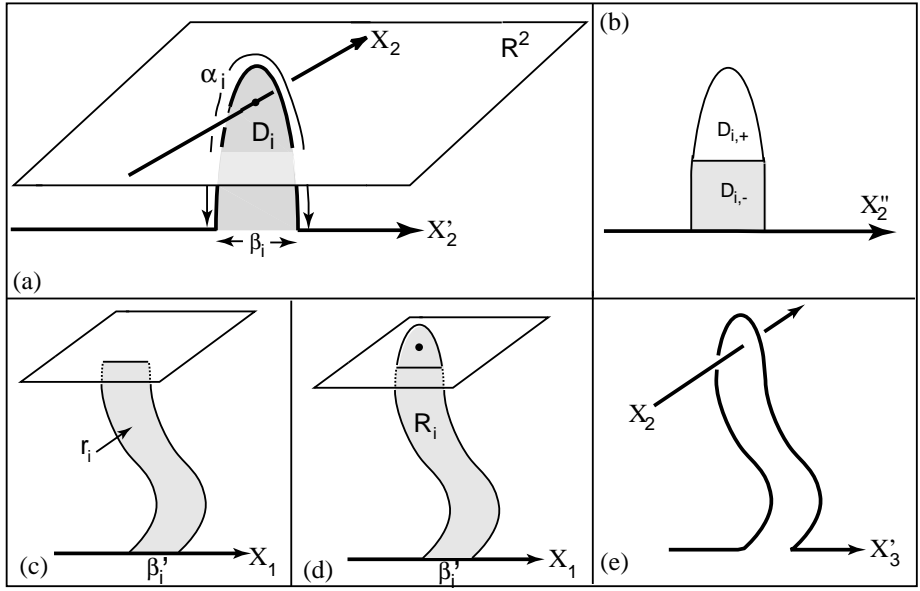 scaled 700}}
\caption{Steps in the construction of $X_3$}
\label{figure:MT3}
\end{figure}

\noindent {\bf Proof of Theorem \ref{theorem:MT}:}   
By hypothesis, we are given closed braids
$X_1$ and $X_2$ which represent the same oriented link type in 3-space.  We must
prove that $X_1\equiv X_2$. We may assume without loss of generality that 
$X_1$ and
$X_2$ are situated in distinct half-spaces (so that they are 
geometrically unlinked),
with
$X_2$ far above $X_1$.  To prove  that $X_1\equiv X_2$ we will 
construct a series of
links $X_2^\prime, X_2^{\prime
\prime}, X_3^\prime, X_3$, where $X_2$ is a preferred longitude for
$X_3$ and $X_3$ is the braid-connected sum of $X_1$ with some number of copies of
the unknot. By Lemma \ref{lemma:MT true for braid connected sum with
unknot} it will follow that $X_1\equiv X_3$. By Lemma 
\ref{lemma:MT true for preferred longitude} it  will follow that $X_3\equiv X_2$.
Therefore we will have proved that $X_1\equiv X_2$.  Figure \ref{figure:MT3} may be
useful to the reader in following the steps of the proof.

To begin the construction choose a Seifert surface ${\bf F}_2$ for $X_2$ and a
preferred longitude $X_2^\prime \subset {\bf F}_2$ for $X_2$.   We will assume
that
$X_2^\prime$ lies in a collar neighborhood of $X_2$, chosen to be small enough so
that $X_2'$ is also a closed braid.  We will also assume that $X_2'$ lies
below $X_2$ everywhere except for little hooks where it is forced to travel
over a strand of $X_2$, as in Figure \ref{figure:MT3}(a), because in general
$X_2\cup X_2'$ is not a split link.  Applying Lemma \ref{lemma:unlinking} we
construct a link
$X_2^\pp$ which is isotopic to $X_2^\prime$ and geometrically unlinked
from $X_2$. The isotopy is a push of $X_2^\prime$ across $k$ disjoint discs
$D_1,\dots,D_k$, replacing each $\alpha_i\subset\partial D_i$ by $\beta_i =
\partial D_i\setminus\alpha_i$. By modifying the subarcs $\beta_i$ of
$X_2^\pp$ a little bit, if necessary, we may assume that $\beta_i$ is
transverse to every fiber
$H_\theta$ of $\cH$, so that $X_2,  X_2'$ and $X_2^\pp$ are all closed braids.

It will be convenient to assume that $X_2$ (resp. $X_2^\pp$ and
$X_1$) lie in the half-spaces $\reals^3_+$ (resp. $\reals^3_-$), where the two
half-spaces are separated by a plane we refer to as $\reals^2$, also that $X_2'$ lies
in $\reals^3_-$ everywhere except for the $k$ hooks where it passes  up and over
$X_2$, intersecting $\reals^2$ twice. See Figure \ref{figure:MT3}(a).  Passing to
Figure \ref{figure:MT3}(b) we think of each disc $D_i, i=1,\dots,k$ as a tall thin
semi-circular disc which is divided by $\reals^2$ into a semi-circular disc
$D_{i,+}\subset\reals^3_+$ and a rectangular disc $D_{i,-}\subset\reals^3_-$,
chosen so that $D_{i,-}\cap X_2^\pp = \beta_i$, where $\beta_i=D_i\cap X_2^\pp$ is
the lower edge of the rectangle $D_{i,-}$.

Noting that $X_2^\pp$ and $X_1$ both represent $\cX$, and that both
are in $\reals^3_-$, we may find a homeomorphism $g:
\reals^3_-\to\reals^3_-$ which is the identity on $\reals^2$ with $g(X_2^\pp)
= X_1$. Extend $g$ by the identity on $\reals^3_+$ to a homeomorphism
$G:\reals^3\to\reals^3$.  Let $r_i = G(D_{i,-})$ and let $R_i = G(D_i) = r_i\cup
D_{i,+}$.  See Figures \ref{figure:MT3}(c) and (d). The facts that (1) $G$ is a
homeomorphism which is the identity in $\reals^3_+$  and (ii) if
$i\not=j$ then $ D_{i,-}\cap D_{j,-} = D_i\cap D_j = \emptyset$ tell us that the
$r_i's$ and  the $R_i's$ are pairwise disjoint embedded discs. By 
construction $r_i\cap X_2 =
\emptyset$, whereas $R_i\cap X_2$ is a single point in the disc $D_{i,+}$.  Each
$r_i$ joins $X_2^\pp$ to $X_1$, meeting $X_2^\pp$ in the arc
$\beta_i$ and $X_1$ in the arc
$\beta_i' = G(\beta_i)$.  Each $R_i$ joins $X_2'$ to $X_1$,
meeting $X_2'$ in $\alpha_i$ and $X_1$ in $\beta_i$.  

Let $X_3^\prime$ be the
link which is obtained from $X_1$ by replacing each $\beta_i'\subset X_1$ by
$\partial R_i - \beta_i'$.  Then $X_3^\prime$ is constructed from $X_1$
by attaching $k$ pairwise disjoint long thin hooks to $X_1$. See Figure
\ref{figure:MT3}(e). There are two important aspects of our construction: 
\begin{enumerate}
\item $X_2\cup X_3'$ has the same link type as $X_2\cup
X_2'$.  For, by construction, the homeomorphism $G^{-1}:\reals^3\to\reals^3$ sends
$X_2\cup X_3'$ to $X_2\cup X_2'$. 
\item $X_3'$ is the connected sum of $X_1$ and $k$ copies of the
unknot, the $i^{th}$ copy being $\partial R_i$.
\end{enumerate}
Now recall that by our initial
construction we had chosen $X_2^\prime$ to be a preferred longitude for $X_2$. 
From (1) it follows that $X_3^\prime$ is also a preferred
longitude for $X_2$. There are 2 cases. 

\underline{Case 1:} If $X_3'$ is isotopic in the complement of
the axis to a closed braid, then let $X_3$ be that closed braid.  By  construction
$X_3$ is the connected sum of $X_1$ and
$k$  closed-braid copies of the unknot.  By Lemma \ref{lemma:MT true for braid
connected sum with unknot} we conclude that $X_1\equiv X_3$. 
We already know that $X_3$ is a preferred longitude for $X_2$. Changing our
point of view, it follows that $X_2$ is a preferred longitude for $X_3$. 
Choose a  Seifert surface
${\bf F}_3$ for $X_3$.  Holding $X_3$ fixed, we may then isotope the interior of
${\bf F}_3$  so that $X_2$  lies on
${\bf F}_3$. But then $X_2$ and $X_3$ cobound a subannulus 
$\cA\subset{\bf F}_3$. But then, by Lemma
\ref{lemma:MT true for preferred longitude}, we conclude that $X_3\equiv X_2$.
Therefore  $X_1\equiv X_2$. 

\begin{figure}[htpb]
\centerline{\BoxedEPSF{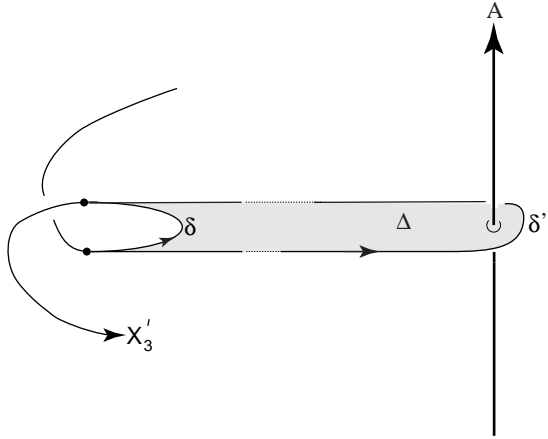 scaled 500}}
\caption{Alexander's trick}
\label{figure:MT4}
\end{figure}

\underline{Case 2:} In general
$X_3^\prime$ will not be a closed braid because the two arcs in
$\partial r_i - \beta_i\cup\beta_i'$ will in general not be transverse to the
fibers of $\cH$.  In this case we  apply Alexander's trick (see
\cite{Alexander}) to change the subarcs that are `in braid position' to a union of arcs
which are in braid position. Alexander does this by subdiving any such
wrongly ordered arc $\delta$ into appropriate smaller arcs, each of which can be
pushed across a disc $\Delta$ which intersects the axis ${\bf A}$ once, as in Figure
\ref{figure:MT4}. 
In this construction Alexander shows that one may always choose the disc $\Delta$ so
that its interior has empty intersection with $X_3'$.  This is accomplished by
choosing the point ${\bf A}\cap\Delta$ to have very large or very small
$z$-coordinate.  In our situation the fact
all of the arcs which are not in braid position are in $\reals^3_-$, whereas $X_2$ is
in $\reals^3_+$ allows us to assume that all such modifications occur far away from
$X_2$. Therefore $X_2\cup X_3$ has the same link type as
$X_2\cup X_3'$.   By the argument for case (1), it then follows that the
proof of Theorem \ref{theorem:MT} is complete. \endpf  

\np{\bf Remark 1:} If $V$ and $W$ are closed braid representatives of knots $\cV$ and
$\cW$, and if $X=V\# W$ represents $\cV\#\cW$, we say that $X$ is the {\it
braid connected sum} of $V$ and $W$ if there is a 2-sphere in 3-space which intersects
$X$ twice and intersects the axis twice, separating the closed braids $V$ and $W$. 
We stated in
$\S$\ref{section:introduction} that our proof of Theorem 1 would show that $X_3$ is
obtained from $X_1$ by taking the braid connected sum of $X_1$ with some number of
closed braid representatives of the unknot $\cU$, say $U_1,\dots,U_k$.  We have,
indeed, shown that there are closed representatives $U_1,\dots,U_k$ of the unknot and
that $X_3=X_1\# U_1\#\cdots\# U_k$, but braid connected sum seems like a
stronger assertion. However, the main result of \cite{BM-4} asserts that in this
situation, after a sequence of exchange moves and isotopies in the complement of the
braid axis, we may assume that our connected sum of closed braids is in fact the braid
connected sum.

\

\np {\bf Remark 2:} In our proof of Theorem \ref{theorem:MT} we only treated the
case of knots, but the proof works equally well for
links. A detailed discussion of the case of links (and of many other aspects of the
construction given here) will appear in the forthcoming manuscript \cite{BM-8}.

\end{document}